\documentclass[11pt,tgrind,psbox]{article}
\usepackage{latexsym}
\usepackage{graphics}
\usepackage{amsmath}
\usepackage{xspace}
\usepackage{amssymb}
\usepackage{psfrag}
\usepackage{epsfig}
\usepackage{fullpage}
\usepackage{amsmath,amsthm}
\usepackage{epsfig}
\usepackage{pst-all}

\newcommand {\bel}[1]{\begin{align*}}
\newcommand {\eel}[1]{\end{align*}}
\newcommand {\bea}{\begin{eqnarray}}
\newcommand {\eea}{\end{eqnarray}}

\newcommand{\Perm}{\ensuremath{\operatorname{Perm}}\xspace}

\newcommand{\G}{\mathbb{G}}

\newcommand{\ignore}[1]{\relax}

\newtheorem{theorem}{Theorem}

\newtheorem{lemma}{Lemma}

\newtheorem{coro}{Corollary}

\newtheorem{Defi}{Definition}

\title{A Deterministic Approximation Algorithm for Computing a Permanent of a $0,1$ matrix}

\author{
 {\sf David Gamarnik }
  \thanks{Operations Research Center and Sloan School of Management, MIT, Cambridge, MA,  02139, e-mail: {\tt
gamarnik@mit.edu}}\and
{\sf Dmitriy Katz} \thanks{Operations Research Center, MIT, Cambridge, MA,  02139, e-mail: {\tt
dimdim@mit.edu}}}

\begin{document}

\maketitle

\begin{abstract}
We construct a deterministic approximation algorithm for computing a permanent of a $0,1$
$n$ by $n$ matrix  to within a multiplicative factor $(1+\epsilon)^n$, for arbitrary $\epsilon>0$.
When the graph underlying the matrix is a constant degree expander our algorithm runs in polynomial time (PTAS).
In the general case the running time of the algorithm is $\exp(O(n^{2\over 3}\log^3n))$.
For the class of graphs which are constant degree expanders
the first result is an improvement over the best known approximation factor $e^n$ obtained in \cite{LinialSamorodnitskyWigderson}.

Our results use a recently developed deterministic approximation  algorithm for counting partial matchings of a
graph \cite{BayatiGamarnikKatzNairTetali} and Jerrum-Vazirani expander decomposition method of \cite{JerrumVazirani}.
\end{abstract}


\section{Introduction}
A permanent of an $n$ by $n$  matrix $A=(a_{i,j})$ is $\Perm(\G)\triangleq\sum_{\sigma}\prod_{1\le i\le n}a_{i,\sigma(i)}$,
where $\sigma$ runs over the elements of the permutation group on the set $1,2,\ldots,n$.
When $A$ is a zero-one matrix, $\Perm(\G)$ counts the number of perfect matching in the graph corresponding to the
adjacency matrix $A$. Permanent
is naturally related to the determinant of  $A$ (signs in front of products are removed). Yet, while
the determinant of a matrix can be computed in polynomial time, the problem of computing a permanent
belongs to the $\#P$ class even when $A$ is a zero-one matrix \cite{valiantPermanent}.
Thus, modulo a basic complexity theoretic conjecture, no polynomial time algorithm
exists for computing permanents.

A lot of research was devoted to constructing approximation algorithms for computing a permanent of a matrix.
The first breakthrough came from Jerrum and Sinclair \cite{JerrumSinclair} who constructed
a randomized fully polynomial time approximation scheme (FPTAS) for $0,1$
matrices satisfying certain ``polynomial growth'' condition.
This condition relates the ratio of perfect to near-perfect matchings and requires that this ratio
is at most a polynomial in $n$. The algorithm is based on a rapidly mixing Markov chain which runs on
perfect and near perfect matchings.
Using this algorithm as a black-box Jerrum and Vazirani constructed a randomized
approximation algorithm for an arbitrary zero-one matrix, but with mildly exponential running time $\exp(O(n^{1\over 2}\log^2 n))$.
A recent dramatic improvement was obtained by Jerrum, Sinclair and Vigoda \cite{JerrumSinclairVigoda},
who constructed FPTAS for an arbitrary matrix with non-negative entries.

Unfortunately, the randomization aspect of the algorithm of \cite{JerrumSinclairVigoda} is quite crucial.
It is not known how to derandomize their algorithm and the best known deterministic approximation algorithm
was given by Linial, Samoridnitsky and Wigderson \cite{LinialSamorodnitskyWigderson}, who provide only
$e^n$ multiplicative approximation factor. Their algorithm is based on a FPTAS for a related matrix scaling problem and uses
van der Waerden's conjecture, which states that a permanent of every doubly stochastic matrix is at least $n!/n^n$.
The $e^n$ approximation factor can be improved to $(k/(k-1))^{kn}$ for the case of of matrices with row and column sums
bounded by $k$, using the Schrijver's bound \cite{SchrijverBound}. For generalizations and simplified proof of Schrijver's
bound see Gurvits~\cite{GurvitsSB}.
Thus approximation of a permanent
is one of the famous algorithmic problems where the existing gap between the randomized and deterministic algorithms is so profound.

In this paper we establish the following two results.
First we construct a polynomial time algorithm which for every $\epsilon$ provides
$(1+\epsilon)^n$ multiplicative approximation factor for a permanent of a $0,1$ matrix, when the underlying graph
is a constant degree expander. The definition
of the latter is given in the subsequent section. While our algorithm is polynomial, it is not fully polynomial, as the term
$1/\epsilon$ appears in the exponent of the running time. Thus we significantly improve the $e^n$ factor
of \cite{LinialSamorodnitskyWigderson} for this class of graphs. Next we construct an  algorithm providing
the same approximation factor $(1+\epsilon)^n$ for a permanent of an \emph{arbitrary} $0,1$ matrix.
The running time of the algorithm
$\exp(O(n^{2\over 3}\log^3n))$. The main technical ingredient of our results is the reduction of the problem to computing
a partition function (see the next section for the definition) corresponding to the collection of all partial and full
matchings of the graph. Recently there was a progress in constructing deterministic approximation algorithms for various
counting problems without the usage of Markov chain based methods. The approach was introduced by
Bandyopadhyay and Gamarnik~\cite{BandyopadhyayGamarnikCounting}, and Weitz~\cite{weitzCounting} and is based on
correlation decay property originating in statistical physics literature. Subsequent works in this direction
include Gamarnik and Katz~\cite{GamarnikKatz} and Bayati et al. \cite{BayatiGamarnikKatzNairTetali}. It is
the second of these two works which provides the basis of our work. Specifically,
a polynomial time approximation  algorithm is constructed in \cite{BayatiGamarnikKatzNairTetali}
for computing the partition function corresponding to partial matchings
of a an arbitrary constant degree graph. The same algorithm solves the underlying problem for an arbitrary graph but
with the running time $\exp(O(n^{1\over 2}\log^2n))$. In this paper obtain appropriate bounds on the permanent
of a matrix in terms of the partition function of partial matchings of a graph. Additionally, our
second result corresponding to arbitrary graphs relies heavily
on the expander decomposition technique developed by Jerrum and Vazirani~\cite{JerrumVazirani} for constructing
a mildly exponential approximation algorithm for computing a permanent of a $0,1$ matrix.

The remainder of the paper is organized as follows. Definitions, preliminaries and the statements of the main results
are subject of the following section. Section~\ref{section:ConstantExpander} is devoted to constructing approximation
algorithm for the case when the underlying graph is a
constant degree expander. The general case is considered in Section~\ref{section:General}. We conclude with some
possible further directions in Section~\ref{section:conclusions}.

\section{Preliminaries and the main result}\label{section:MainResult}
Consider a simple undirected $n$ by $n$ bi-partite graph  $\G$ with the node
set $V=V_1\cup V_2$, $|V_1|=|V_2|=n$.
Let $E$ be
the set of  edges of the graph. $N(v,\G)\subset V$ denotes the set of neighbors of $v\in V$
and $\Delta(v)=|N(v,\G)|$ denotes the degree
of the node $v$. The degree of the graph is $\Delta\triangleq \max_v \Delta(v)$.
Given a set of nodes $A$, we denote by $N(A)$ or specifically by $N(A,\G)$ the set of neighbors of $A$ (in $\G$).
Given $\alpha>0$, a graph $\G$ is defined to be an \emph{expander} or specifically $\alpha$-\emph{expander} if
for every set of nodes $A\subset V_i, i=1,2$ such that $|A|\leq n/2$, the inequality $|N(A)|\geq (1+\alpha)|A|$ holds.
We also define
\begin{align*}
\alpha(\G)\triangleq \max_{A}{|N(A)|\over |A|}-1
\end{align*}
to be the \emph{expansion} of the graph $\G$, where the maximum is over all subsets $A\subset V_i, i=1,2$
with $|A|\le n/2$. Clearly, $\G$ is $\alpha$-expander if its expansion is at least $\alpha$.

A matching is a subset $M\subset E$ such that no two edges in $M$ share a node.
For every $k\leq n$ let $M(k)$ be the number of size $k$ matchings in $\G$.
Specifically, $M(n)$ is the number of full matchings.

Given a graph $\G$ let $A=(a_{ij})$ be the corresponding adjacency matrix. The rows and columns of $A$
are indexed by nodes of $V_1$ and $V_2$ respectively, and $a_{ij}=1$ if $(v_i,v_j)\in E$
and $a_{ij}=0$ otherwise. It is immediate that $M(n)=\Perm(A)$.

A parameter $\lambda>1$ is fixed called the \emph{activity}. The \emph{partition function} corresponding to $\lambda$
is defined as
\begin{align*}
Z(\lambda,\G)=\sum_{M}\lambda^{|M|}=\sum_{1\le k\le n}\lambda^k M(k).
\end{align*}
A partition function is an important object in statistical physics. The case of matching is usually called
monomer-dimer model in the statistical physics literature.
\begin{Defi}
An approximation algorithm ${\cal A}$ is defined to be a Fully Polynomial Time Approximation Scheme (FPTAS) for a computing
$Z(\lambda,\G)$ if given arbitrary $\delta>0$ it produces a value $\hat Z$ satisfying
\begin{align*}
{1\over 1+\delta}\leq {\hat Z\over Z(\lambda,\G)}\leq 1+\delta,
\end{align*}
in time which is polynomial in $n$ and ${1\over \delta}$.
\end{Defi}

The following result was established in \cite{BayatiGamarnikKatzNairTetali}.
\begin{theorem}\label{theorem:Matchings}
There exist a deterministic algorithm which provides a FPTAS for computing $Z(\lambda,\G)$ for an arbitrary
graph/activity  pair $\G,\lambda$ when $\Delta$ and $\lambda$ are constants. The complexity of the same
algorithm is $\exp(O(\sqrt{\lambda n}\log^2n))$ for general $\Delta$ and $\lambda$.
\end{theorem}
The case $\lambda=1$ corresponds to counting the number of partial matchings of a graph.
In this paper we use the algorithm underlying Theorem~\ref{theorem:Matchings} as a subroutine to devise
an approximation algorithm for computing a permanent. For this purpose we will be making $\lambda$ to be
appropriately large. Throughout the paper we assume $\lambda>10$.


We now state the main two results of this paper. In the next and the following results
the notion of  $(1+ \epsilon)^n$ multiplicative approximation factor of $\Perm(\G)$ corresponds to obtaining
a value $\hat Z$ satisfying $(1+\epsilon)^{-n}\le {\hat Z\over \Perm(\G)}\le (1+\epsilon)^n$.
\begin{theorem}\label{theorem:MainResultPermConst}
Let $\G$ be an $n$ by $n$ bi-partite $\alpha$-expander and let $\epsilon>0$.
There
exist a deterministic  $(1+ \epsilon)^n$ approximation algorithm for computing $\Perm(\G)$ with
running time  $\exp(O(\sqrt{\epsilon^{-1}\alpha^{-1}n}\log^3n))$. Moreover, the running time is polynomial in $n$
whenever $\Delta, \alpha$ are constants.
\end{theorem}
Our algorithm corresponding to the second part of the theorem,
while polynomial, is not strongly polynomial. As we shall see, the dependence of the running time
on the approximation parameter $\epsilon$ is of the form $n^{O({1\over \epsilon})}$. While our approximation factor
$(1+ \epsilon)^n$ is a far cry from PTAS (namely $1+\epsilon$ approximation factor),
it is still a significant improvement over $e^n$ factor constructed
in \cite{LinialSamorodnitskyWigderson}.

Our second result does not require any restrictions on the underlying graph.
\begin{theorem}\label{theorem:MainResultPermGen}
There exist a deterministic  $(1+ \epsilon)^n$ approximation algorithm for computing $\Perm(\G)$ of an arbitrary
$n$ by $n$ bi-partite  graph $\G$ which runs in time $\exp(O(\epsilon^{-{1\over 2}}n^{2\over 3}\log^3n))$.
\end{theorem}
Thus, similarly to \cite{JerrumVazirani}, our algorithm provides a mildly exponential algorithm for approximating
a permanent (with a weaker approximation factor $(1+ \epsilon)^n$).

\section{Constant degree expanders}\label{section:ConstantExpander}
Proof of Theorems~\ref{theorem:MainResultPermConst} is given in this section.
We begin by establishing some preliminary results. We assume without the loss
of generality that $M(n)\ge 1$.
Consider an arbitrary $k$ matching $M$ between sets $A_1\subset V_1, A_2\subset V_2, |A_1|=|A_2|=k$.
A path $v_1,v_2,\ldots,v_{2r}$ is defined to be an \emph{alternating path} wrt $M$ if
$v_1\in V_1\setminus A_1,v_{2r}\in V_2\setminus A_2$, if $(v_{2l},v_{2l+1})\in M, 1\le l\le r-1$ and
$(v_{2l-1},v_{2l})\notin M, 1\le l\le r$. Observe that given $M$ and an alternating path $P$
one can construct a $k+1$ matching,  by subtracting from $M$ all the edges
$(v_{2l},v_{2l+1}), 1\le l\le r$ and adding all the edges $(v_{2l-1},v_{2l}), 1\le l\le r-1$.
The length of this alternating path is defined to be $2r$.

\begin{lemma}\label{lemma:alternatinglength}
Let $M$ be a $k$ matching between $A_1$ and $A_2$ for $k\le n-1$.
For every  set $L\subset V_1\setminus A_1$ with $|L|\ge (n-k)/2$
there exists an alternating path $P$ with end points in $L$ and $V_2\setminus A_2$ and with length at most
\begin{align}\label{eq:alternatinglength}
O\Big(\log({n\over n-k})\log^{-1}(1+\alpha )\Big).
\end{align}
\end{lemma}

\begin{proof}
The proof is similar to  the argument of Lemma 2 \cite{JerrumVazirani}.
Let $R=V_2\setminus A_2$.
Let $L_r$ be the set of nodes in $V_2$ reachable from $L$ en route of alternating paths with length
at most $2r$. Let $r_0$ be a minimum $r$ satisfying $(1+\alpha)^r(n-k)/2>n$.
Note $r_0=O\Big(\log({n\over n-k})\log^{-1}(1+\alpha )\Big)$.
Since $|L|\ge (n-k)/2$,
then $l(r_0)\triangleq\min((1+\alpha )^{r_0}|L|,{n\over 2}+1)={n\over 2}+1$.
For every $r\le r_0$, by the
expansion property either $|L_r|\geq \min((1+\alpha )^r|L|,{n\over 2}+1)$, or
$L_{r'}\cap R\neq \emptyset$ for some $r'\le r$. In the second case we found an alternating path with
length $\le l(r_0)$. In the first case we have $|L_{r_0}|>n/2$. We now claim that for every $r>r_0$
either $L_r\cap V_2\setminus A_2\neq\emptyset$ or $|V_2\setminus L_r|\le {1\over (1+\alpha)^{r-r_0}}(n/2)$.
Indeed, let $L_r'\subset V_1$ be the set of nodes matched to $L_r$. In particular $|L_{r'}|>n/2$.
Since $|L_{r+1}|=|N(L_{r'})|>n/2$, then
\begin{align*}
|V_2\setminus L_r|=|V_1\setminus L_{r'}|
\ge |N(V_2\setminus L_{r+1})|\geq (1+\alpha)|V_2\setminus L_{r+1}|,
\end{align*}
and the assertion is established by induction. For $r\ge 2r_0$ we obtain
${1\over (1+\alpha)^{r-r_0}}(n/2)<n-k$ and thus we must have $L_r\cap V_2\setminus A_2\neq \emptyset$.
We established that for some $r\le 2r_0$ there exists an alternating path with an end points
in  $L$ and $V_2\setminus A_2$. The length of this path is $O(r_0)=O({n\over n-k}\log^{-1}(1+\alpha ))$.
This completes the proof.
\end{proof}

\begin{lemma}\label{lemma:Mkn}
For every $k\le n-1$
\begin{align}\label{eq:Mkk+1}
{M(k)\over M(k+1)}\leq 2\Big({n\over n-k}\Big)^{O\Big(\log^{-1}(1+\alpha )\log\Delta \Big)}.
\end{align}
As a result
\begin{align}\label{eq:Mkn}
{M(k)\over M(n)}\leq \Big({2en\over n-k}\Big)^{O\Big((n-k)\log^{-1}(1+\alpha )\log\Delta\Big)}.
\end{align}
\end{lemma}

\begin{proof}
Fix an arbitrary $k$ matching $M$ between $A_1\subset V_1$ and $A_2\subset V_2$.
We claim that there exist at least $(n-k)/2$ ~~$k+1$-matchings obtained
from $M$ via an alternating path with length at most the value given by (\ref{eq:alternatinglength}).
Indeed, consider the set of all nodes $v$ in $V_1\setminus A_1$ such that the shortest alternating path
starting from $v$ is larger than the required bound.
By Lemma~\ref{lemma:alternatinglength} this set contains less than $(n-k)/2$ nodes and the assertion follows.

Now consider the following bi-partite graph. The nodes on the left (right) are all $k$ ($k+1$)-matchings.
We put an edge between two matchings $M,M'$ if $M'$ is obtained from $M$ via an alternating path with length
at most (\ref{eq:alternatinglength}). The total number of edges in this graphs is at least  $M(k)(n-k)/2$ by our
observation above. For every matching $M'$ on the right side of the graph the number of edges pointing to it
is at most the number of alternating paths with length at most (\ref{eq:alternatinglength}) which result in
$M'$. For every possible starting node of an alternating path, the number of such alternating paths is at most
\begin{align*}
\Delta^{O\Big(\log({n\over n-k})\log^{-1}(1+\alpha )\Big)}.
\end{align*}
The number of starting nodes is bounded by $n$. Then the total number of edges in this graph is
at most $M(k+1)n\Delta^{O\Big(\log({n\over n-k})\log^{-1}(1+\alpha )\Big)}$. We conclude that
\begin{align*}
M(k+1)n\Delta^{O\Big(\log({n\over n-k})\log^{-1}(1+\alpha )\Big)}\geq M(k)(n-k)/2.
\end{align*}
The  bounds (\ref{eq:Mkk+1}) then follows. From this bound we also obtain
\begin{align*}
{M(k)\over M(n)}\leq 2^{n-k}\Big({n^{n-k}\over (n-k)!}\Big)^{O\Big(\log^{-1}(1+\alpha )\log\Delta \Big)}
=\Big({(2e n)^{n-k}\over (n-k)^{n-k}}\Big)^{O\Big(\log^{-1}(1+\alpha )\log\Delta \Big)},
\end{align*}
where the Stirling's approximation was used in the equality. This is the required bound (\ref{eq:Mkn}).
\end{proof}

\begin{coro}\label{coro:nlambda}
The following holds
\begin{align}\label{eq:twosidebounds}
1\le {Z(\lambda,\G)\over \lambda^n\Perm(\G)}\le \exp(O(n\lambda^{-1}\log^{-1}(1+\alpha )\log\Delta)).
\end{align}
\end{coro}


\begin{proof}
The inequality $Z(\lambda,\G)\ge \lambda^n\Perm(\G)$ is immediate. We focus on the second inequality in
(\ref{eq:twosidebounds}). We need to analyze the ratio
\begin{align}\label{eq:Mratio}
{\lambda^{k}M(k)\over \lambda^n M(n)}=\lambda^{-(n-k)}{M(k)\over M(n)}.
\end{align}
We set $c_n\triangleq \log^{-1}(1+\alpha )\log\Delta$ for simplicity. Applying the second part of Lemma~\ref{lemma:Mkn}
\begin{align*}
\lambda^{-(n-k)}{M(k)\over M(n)}\le \exp\Big(O\big(c_n(n-k)\big(\log n-\log(n-k)
-\log{\lambda\over 2e}\big)\big)\Big)
\end{align*}
Consider the problem of maximizing
\begin{align*}
g(x)\triangleq x\log n-x\log x-x\log{\lambda\over 2e}
\end{align*}
in the range $x\in [1,n]$. The boundary cases $x=1,n$ give respectively values
$\log(2en/\lambda),-n\log(\lambda/2e)$. The second quantity is negative when $\lambda>2e$
(recall our assumption $\lambda>10$).
To find another candidate for the largest value, we take the derivative with respect to $x$ and equating it to zero
we obtain
\begin{align*}
\log n-\log x-\log{\lambda\over 2}=0,
\end{align*}
giving $x=2n/\lambda$. Evaluating $g(x)$ at this value simplifies to $2n/\lambda$
and this gives the largest value of $g$ when $n$ is larger than some $\lambda$ dependent constant.
We conclude that the left-hand side of (\ref{eq:Mratio}) is at most $\exp(O({c_n n\over \lambda}))$, implying
\begin{align*}
{Z(\lambda,\G)\over \lambda^n\Perm(\G)}\le
n\exp(O({c_n n\over \lambda}))=\exp(O({c_nn\over \lambda}+\log n))=\exp(O({c_nn\over \lambda})).
\end{align*}

\end{proof}

\begin{proof}[Proof of Theorem~\ref{theorem:MainResultPermConst}]
Fix an arbitrary constant $\epsilon>0$.
We select the smallest  $\lambda$
so that $\exp(O(\lambda^{-1}\log^{-1}(1+\alpha)\log\Delta))<1+\epsilon$.
It is clear that $\lambda=O(\epsilon^{-1}\alpha^{-1}\log\Delta)$.
We compute
an $\epsilon$-approximation $\tilde Z$ of $Z(\lambda,\G)$ using an algorithm from Theorem~\ref{theorem:Matchings}
for computing $Z(\lambda,\G)$.
By Corollary~\ref{coro:nlambda}, it satisfies
\begin{align*}
(1-\epsilon)\le {\tilde Z\over \lambda^n\Perm(\G)}\le (1+\epsilon)^{n+1}.
\end{align*}
Then $\tilde Z/\lambda^n$ provides the required approximation.
The complexity of this algorithm is
\begin{align*}
\exp(O(\sqrt{\lambda n}\log^2n))=\exp(O(\sqrt{\epsilon^{-1}\alpha^{-1}(\log\Delta) n}\log^2n))=
\exp(O(\sqrt{\epsilon^{-1}\alpha^{-1}n}\log^3n)),
\end{align*}
and the first part of the theorem is established.

For the second part we observe that $\lambda=O(\epsilon^{-1}\alpha^{-1}\log\Delta)$ is a constant whenever
$\alpha$ and $\Delta$ are constants. We
recall from Theorem~\ref{theorem:Matchings} that the algorithm
for computing $Z(\lambda,\G)$ is polynomial time, under these assumptions.
\end{proof}

\section{General graphs}\label{section:General}
\begin{proof}[Proof of Theorem~\ref{theorem:MainResultPermGen}]
Our approach borrows heavily from the Jerrum-Vazirani expander decomposition approach \cite{JerrumVazirani}.
The idea is to decompose the underlying graph into a collection of subgraphs with a suitable expansion properties
and apply an algorithm for computing a permanent recursively. In \cite{JerrumVazirani}
the subroutine used is based on the algorithm relying on rapidly mixing Markov chain. Here we use the deterministic
algorithm constructed in the proof of Theorem~\ref{theorem:MainResultPermConst}.

The following result is established in \cite{JerrumVazirani} (Lemma 4).
There exists and algorithm, called {\tt TestExpansion} which on input $\G,\alpha$ either
identifies that $\G$ is an $\alpha$-expander, or identifies a set $A\subset V_1, |A|\le n/2$ such that
$N(A)\le (1+2\alpha)|A|$. The running time of the algorithm is $\exp(O(\alpha n\log n))$.

Given an arbitrary set $A\subset V_1$ it is straightforward to observe that
\begin{align}\label{eq:product}
\Perm(\G)=\sum_{B\subset N(A), |B|=|A|}\Perm(A,B)\Perm(A^c,B^c)
\end{align}
where $A^c=V_1\setminus A,B^c=V_2\setminus B$ and $\Perm(A,B)$ is the permanent of the subgraph induced by
$A$ and $B$. Based on this observation we propose the following recursive algorithm for computing $\Perm(\G)$.
We suppose that we have an algorithm $\mathcal{A}_r$ which computes $(1+ \epsilon)^r$ factor approximation of a permanent
of any $r$ by $r$ bipartite graph in time $g(r)$, for every $r\le n-1$. We use it to construct $\mathcal{A}_n$ as follows.
Set $\alpha =n^{-1/3}$. Run algorithm {\tt TestExpansion} on $\G$. The running time is
$\exp(O(n^{2\over 3}\log n))$.
If the algorithm returns no set $A$, then the
underlying graph is an $\alpha $-expander and we use algorithm of Theorem~\ref{theorem:MainResultPermConst}
to obtain an $(1+ \epsilon)^n$ approximation of $\Perm(\G)$.  The running time
is $\exp(O(\sqrt{\epsilon^{-1}n^{1\over 3}n}\log^3n))=\exp(O(\epsilon^{-{1\over 2}}n^{2\over 3}\log^3n))$, and the overall
running time is $\exp(O(n^{2\over 3}\log n))+\exp(O(\epsilon^{-{1\over 2}}n^{2\over 3}\log^3n))=
\exp(O(\epsilon^{-{1\over 2}}n^{2\over 3}\log^3n))$.
Let $c_0$ be the constant hidden in $O(\cdot)$. From now on, treating
$\epsilon$ as constant, we hide $\epsilon^{-{1\over 2}}$ in the $O(\cdot)$ term. In the end
we observe that the dependence of the running time on $\epsilon$ is of the form $\exp(O(\epsilon^{-{1\over 2}}))$.

Otherwise the algorithm  {\tt TestExpansion} identifies a set
$A$ with $|N(A)|\le (1+2\alpha)|A|$. We estimate $\Perm(A,B)$ and $\Perm(A^c,B^c)$ using algorithm
$\mathcal{A}_r$ with $r=|A|, |A^c|$ respectively.
Then we estimate $\Perm(\G)$ using the expression (\ref{eq:product}).
For every product $\Perm(A,B)\Perm(A^c,B^c)$ our approximation factor
is $(1+\epsilon)^{|A|}(1+\epsilon)^{|A^c|}=(1+\epsilon)^n$ by the recursive assumption.

Now we obtain an upper bound on $g(n)$ and specifically show that it is $\exp(O(n^{2\over 3}\log^3 n))$.
Let $c>c_0$ be a very large constant. By the recursive assumption $g(r)\le \exp(cr^{2\over 3}\log^3 r)), ~r\le n-1$.
The function $g$ satisfies the following bound
\begin{align}\label{eq:gn}
g(n)\le \max\Big(\exp(c_0 n^{2\over 3}\log^3n),\max_{1\le r\le n/2}{r(1+2\alpha ) \choose r}(g(r)+g(n-r))\Big).
\end{align}
Here the term ${r(1+2\alpha ) \choose r}$ comes from performing the computation over all subsets $B$ of $N(A)$ with
size $|B|=|A|$. This term is trivially upper bounded by
\begin{align}
n^{2\alpha  r}&=\exp(2n^{-{1\over 3}}r\log n) \label{eq:alphar}\\
&\le \exp(n^{{2\over 3}}\log n).
\end{align}
Now, by the recursive assumption we have for $r\le n/2$ that
$g(r)\le \exp(cr^{2\over 3}\log^3 r))\le \exp((3/4)cn^{2\over 3}\log^3 n))$, which gives
\begin{align*}
{r(1+2\alpha ) \choose r}g(r)\le \exp((7/8)cn^{2\over 3}\log^3 n)),
\end{align*}
provided that $1+(3/4)c<(7/8)c$. As for $g(n-r)$, we have for sufficiently large $n$
\begin{align*}
{r(1+2\alpha ) \choose r}g(n-r)&\le \exp(2n^{-{1\over 3}}r\log n)\exp(c(n-r)^{2\over 3}\log^3 (n-r))\\
&\le \exp(2n^{-{1\over 3}}r\log n)\exp(c(n-r)^{2\over 3}\log^3 n)\\
&\stackrel{(a)}{\le}  \exp(2n^{-{1\over 3}}r\log n)\exp(cn^{2\over 3}\log^3 n-(cr/2)n^{-{1\over 3}}\log^3n) \\
&\stackrel{(b)}{\le}   \exp(cn^{2\over 3}\log^3 n-2n^{-{1\over 3}})\\
&\le (1-n^{-{1\over 3}})\exp(cn^{2\over 3}\log^3 n),
\end{align*}
where in $(a)$ we use $(n-r)^{2\over 3}\le n^{2\over 3}-{r\over 2} n^{-{1\over 3}}$ (obtained for example
using Taylor's expansion around $n^{2\over 3}$); and in $(b)$ we use $r\ge 1$ and $c>8$.
Using $\exp((7/8)cn^{2\over 3}\log^3 n))\le n^{-{1\over 3}}\exp(cn^{2\over 3}\log^3 n)$
for large $n$, we obtain that for every $1\le r\le n-1$
\begin{align*}
{r(1+2\alpha ) \choose r}(g(r)+g(n-r))\le \exp(cn^{2\over 3}\log^3 n).
\end{align*}
Combining with (\ref{eq:gn}) we obtain the required bound $g(n)\le \exp(cn^{2\over 3}\log^3 n)$.
\end{proof}

\section{Conclusions}\label{section:conclusions}
We proposed a new deterministic approximation algorithm for computing a permanent of an
$n$ by $n$ $0,1$ matrix.
Our algorithm provides a multiplicative approximation factor $(1+\epsilon)^n$ and runs in polynomial time
for a matrix corresponding to a constant degree expander, and in time $\exp(O(n^{2\over 3}\log^3))$ for an arbitrary matrix.
Our algorithm is based on a recent deterministic approximation algorithm for counting the number
of partial matchings of a graph \cite{BayatiGamarnikKatzNairTetali}.
It is natural to try to extend our results in several directions. One possibility is lifting $0,1$ requirement.
This entails extending the result of \cite{BayatiGamarnikKatzNairTetali} to the case of weighted graphs. It is
reasonable to expect that for some special cases such a program will work.
In this case one can try to extend the results of the present paper to the case of general matrix entries.

A second interesting direction is utilizing the ideas in \cite{JerrumSinclairVigoda}, developed in the Markov chain
setting, as a possibility of getting stronger deterministic approximation algorithms for computing a permanent. The technique
developed in \cite{JerrumSinclairVigoda} allow one to deal with the case when the ratio $M(n-1)/M(n)$ is exponentially
large. Since obtaining amenable bounds on the ratios $M(k)/M(n)$ was required for our algorithms to work, this might
be indeed a fruitful direction for further research.

\section*{Acknowledgements}
The authors gratefully acknowledge several important references obtained from
Alexander Barvinok and Leonid Gurvits.

\bibliographystyle{amsalpha}

\newcommand{\etalchar}[1]{$^{#1}$}
\providecommand{\bysame}{\leavevmode\hbox to3em{\hrulefill}\thinspace}
\providecommand{\MR}{\relax\ifhmode\unskip\space\fi MR }
\providecommand{\MRhref}[2]{%
  \href{http://www.ams.org/mathscinet-getitem?mr=#1}{#2}
}
\providecommand{\href}[2]{#2}

\end{document}